# Burns' equivariant Tamagawa invariant $T\Omega^{loc}(N/\mathbf{Q}, 1)$ for some quaternion fields

Victor Snaith

## 1 Introduction

Inspired by the work of Bloch and Kato in [2], David Burns constructed several "equivariant Tamagawa invariants" associated to motives of number fields. These invariants lie in relative $K$-groups of group-rings of Galois groups and in [3] Burns made several conjectures (see Conjecture 2.2) about their values. In this paper I shall verify Burns' conjecture concerning the invariant $T\Omega^{loc}(N/\mathbf{Q}, 1)$ for some families of quaternion extensions $N/\mathbf{Q}$. The family of quaternion fields which are not covered here will be covered in [8].

The paper is arranged in the following manner. In §2 the relative K-groups of $\mathbf{Z}[G]$ are introduced together with the principle method for constructing elements in them and the Hom-description for representing elements in the case when $G = Q_8$, the quaternion group of order eight. In §3 the formula for $T\Omega^{loc}(N/\mathbf{Q}, 1)$ is given together with several simplifying observations which are special to quaternion extensions. This formula involves a sum of "local terms" whose Hom-description representatives are calculated in §4. In §5, when $N/\mathbf{Q}$ is a quaternion extension whose decomposition group at 2 is a proper subgroup of $Q_8$, the Hom-descriptions of §3 and §4 are combined to prove Conjecture 2.2, which asserts that $T\Omega^{loc}(N/\mathbf{Q}, 1)$ vanishes.

I am very grateful to David Burns for telling me of his work [3]. This paper was written in 1999, during a visit by Werner Bley to the University of Southampton. This one-month visit was sponsored by the EU TMR grant entitled "Algebraic K-theory and linear algebraic groups". I am particularly indebted to Werner Bley for patiently explaining to me the (rather intricate) definition of $T\Omega^{loc}(L/K, 1)$ without which I could not have made these calculations.

## 2 Relative K-groups of group-rings

**2.1** *The exact sequence*

Let $G$ be a finite group. Let $R$ be an integral domain and let $f : \mathbf{Z}[G] \longrightarrow R[G]$ be the corresponding inclusion of group-rings. There is an associated exact sequence of algebraic K-groups, constructed in ([17] p.216), of the form

$$\ldots \longrightarrow K_1(\mathbf{Z}[G]) \longrightarrow K_1(R[G]) \xrightarrow{\delta_1} K_0(\mathbf{Z}[G], R) \xrightarrow{\pi_*} K_0(\mathbf{Z}[G]) \xrightarrow{f_*} K_0(R[G]).$$

Here we have adopted the notation of [3] for the relative K-group in the middle, which was originally denoted by $K_0(\mathbf{Z}[G], f)$ in [17].

An arbitrary element of $K_0(\mathbf{Z}[G], R)$ is represented by a triple $[A, \phi, B]$ in which $A$ and $B$ are finitely generated projective (left) $\mathbf{Z}[G]$-modules and $\phi$ is an isomorphism of $R[G]$-modules of the form

$$\phi : A \otimes R \xrightarrow{\cong} B \otimes R.$$

A presentation for $K_0(\mathbf{Z}[G], R)$ in terms of these triples is given in ([17] p.215). The homomorphisms in the sequence are defined by $\pi_*[A, \phi, B] = [A] - [B] \in K_0(\mathbf{Z}[G])$, $f_*[A] = [A \otimes R] \in K_0(R[G])$ and if $X \in GL_n(R[G])$ is an invertible matrix representing $x \in GL(R[G])^{ab} = K_1(R[G])$ then $\delta_1(x) = [R[G]^n, X, R[G]^n] \in K_0(\mathbf{Z}[G], R)$.

Now consider the situation when $R = \mathbf{Q}$, the rational field, which is the case of particular interest to us. Incidentally, in this case the relative K-theory exact sequence is just the low-dimensional end of a localisation sequence (cf. [11] §5 Theorem 5). The categorical construction of the localisation sequence in [11] shows that there is an isomorphism of the relative K-groups of the form

$$K_0(\mathbf{Z}[G]; \mathbf{Q}) \cong \oplus_{p \ prime} K_0(\mathbf{Z}_p[G]; \mathbf{Q}_p)$$

sending $[A, \phi, B]$ to $\{[A \otimes \mathbf{Z}_p, \phi, B \otimes \mathbf{Z}_p] \mid p \ prime\}$. Here $K_0(\mathbf{Z}_p[G]; \mathbf{Q}_p)$ is defined in an analogous manner to $K_0(\mathbf{Z}[G]; \mathbf{Q})$, replacing $\mathbf{Z}$ by $\mathbf{Z}_p$ and $\mathbf{Q}$ by $\mathbf{Q}_p$.

**2.2** *Constructing elements of $K_0(\mathbf{Z}[G], \mathbf{Q})$*

We shall require methods by which to construct and compare elements in $K_0(\mathbf{Z}[G], \mathbf{Q})$.

Suppose that we have a bounded chain complex of perfect (i.e. finitely generated, projective) $\mathbf{Z}[G]$-modules of the form

$$P^* : \ 0 \longrightarrow P^n \xrightarrow{d_n} P^{n+1} \xrightarrow{d_{n+1}} \ldots \xrightarrow{d_{m-1}} P^m \longrightarrow 0$$

together with a given isomorphism

$$\psi : \oplus_j \ H^{2j+1}(P^*) \otimes \mathbf{Q} \xrightarrow{\cong} \oplus_j \ H^{2j}(P^*) \otimes \mathbf{Q}.$$

Associated to this data is a well-defined element of $K_0(\mathbf{Z}[G], \mathbf{Q})$, $[P^{od}, \phi, P^{ev}]$ in the following manner. Set $P^{od} = \oplus_j P^{2j+1}$ and $P^{ev} = \oplus_j P^{2j}$.

Choose $\mathbf{Q}[G]$-modules splittings for the right-hand epimorphisms in each of the short exact sequences:

$$0 \longrightarrow Ker(d_j) \otimes \mathbf{Q} \longrightarrow P^j \otimes \mathbf{Q} \longrightarrow d_j(P^j) \otimes \mathbf{Q} \longrightarrow 0,$$

$$0 \longrightarrow d_{j-1}(P^{j-1}) \otimes \mathbf{Q} \longrightarrow Ker(d_j) \otimes \mathbf{Q} \longrightarrow H^j(P^*) \otimes \mathbf{Q} \longrightarrow 0.$$

Using these chosen splittings the rational isomorphism $\phi$ is given by the following composition

$$\oplus_j \ P^{2j+1} \otimes \mathbf{Q}$$

$$\cong (\oplus_j \ Ker(d_{2j+1}) \otimes \mathbf{Q}) \oplus (\oplus_j \ d_{2j+1}(P^{2j+1}) \otimes \mathbf{Q})$$

$$\cong (\oplus_j \ d_{2j}(P^{2j}) \otimes \mathbf{Q}) \oplus (\oplus_j \ H^{2j+1}(P^*) \otimes \mathbf{Q}) \oplus (\oplus_j \ d_{2j+1}(P^{2j+1}) \otimes \mathbf{Q})$$

$$\stackrel{\psi}{\cong} (\oplus_j \ d_{2j}(P^{2j}) \otimes \mathbf{Q}) \oplus (\oplus_j \ H^{2j}(P^*) \otimes \mathbf{Q}) \oplus (\oplus_j \ d_{2j+1}(P^{2j+1}) \otimes \mathbf{Q})$$

$$\cong (\oplus_j \ d_{2j}(P^{2j}) \otimes \mathbf{Q}) \oplus (\oplus_j \ Ker(d_{2j}) \otimes \mathbf{Q})$$

$$\cong \oplus_j \ P^{2j} \otimes \mathbf{Q}.$$

The Euler characteristic of the perfect complex is defined to be $\chi(P^*) = \sum_j (-1)^{j+1}[P^j] \in K_0(\mathbf{Z}[G])$ so that $\pi_*[P^{od}, \phi, P^{ev}] = \chi(P^*)$.

**2.3** *The Hom-description*

The Hom-description enables us to represent elements of $K_0(\mathbf{Z}[G], \mathbf{Q})$ in terms of functions on the representation ring of $G$, $R(G)$, given by the Grothendieck group of finite-dimensional, complex representations of $G$. Let $\Omega_{\mathbf{Q}_p}$ denote the absolute Galois group of the $p$-adic rationals $\mathbf{Q}_p$ and let $\overline{\mathbf{Q}}_p$ be an algebraic closure. Choosing an isomorphism of $R(G)$ with the Grothendieck ring of $\overline{\mathbf{Q}}_p$-representations we can endow $R(G)$ with an action by $\Omega_{\mathbf{Q}_p}$. Also the multiplicative group, $(\overline{\mathbf{Q}}_p)^*$, has a natural action by $\Omega_{\mathbf{Q}_p}$. We have a Hom-description isomorphism of the form

$$K_0(\mathbf{Z}_p[G]; \mathbf{Q}_p) \cong \frac{Hom_{\Omega_{\mathbf{Q}_p}}(R(G), (\overline{\mathbf{Q}}_p)^*)}{Det(\mathbf{Z}_p[G]^*)}$$

which is defined in the following manner. Every element of $K_0(\mathbf{Z}_p[G]; \mathbf{Q}_p)$ may be written as $[A, \phi, B]$ in which $A, B$ are finitely generated, free $\mathbf{Z}_p[G]$-modules. We may choose $\mathbf{Z}_p[G]$-bases $\underline{x}_p = \{x_1, x_2, \ldots, x_n\}$ and $\underline{y}_p = \{y_1, y_2, \ldots, y_n\}$ for $A$ and $B$ respectively. Now let $\Phi_p \in GL_n(\mathbf{Q}_p[G])$ denote the matrix which represents the isomorphism, $\phi : A \otimes \mathbf{Q}_p \xrightarrow{\cong} B \otimes \mathbf{Q}_p$. Hence $\phi(\underline{x}_p) = \Phi_p \underline{y}_p$. Now suppose that $\rho : G \longrightarrow GL_m(\overline{\mathbf{Q}}_p)$ is a representation of $G$. Define $Det(\Phi_p)(\rho) \in (\overline{\mathbf{Q}}_p)^*$ by the formula

$$Det(\Phi_p)(\rho) = det(\rho(\Phi_p))$$

where $\rho(\Phi_p) \in GL_{nm}(\overline{\mathbf{Q}}_p)$ is obtained by applying $\rho$ to the group elements in $\Phi_p$. This formula uniquely characterises a Galois equivariant homomorphism

$$Det(\Phi_p) \in Hom_{\Omega_{\mathbf{Q}_p}}(R(G), (\overline{\mathbf{Q}}_p)^*)$$

which gives a Hom-description representative of $[A, \phi, B] \in K_0(\mathbf{Z}_p[G]; \mathbf{Q}_p)$.

Performing this construction prime by prime gives rise to a Hom-description of $K_0(\mathbf{Z}[G]; \mathbf{Q})$ in term of Galois equivariant homomorphisms from $R(G)$ to $J^*(\overline{\mathbf{Q}})$, the idèles of $\overline{\mathbf{Q}}$, an algebraic closure of the rationals. The class-group, $\mathcal{CL}(\mathbf{Z}[G])$, is defined to be the reduced $K_0$-group

$$\mathcal{CL}(\mathbf{Z}[G]) = ker(rank : K_0(\mathbf{Z}[G]) \longrightarrow \mathbf{Z}).$$

A similar construction with local bases yields a Hom-description isomorphism of the form ([6], [14] p.40, [15] p.115)

$$Det : \mathcal{CL}(\mathbf{Z}[G]) \xrightarrow{\cong} \frac{Hom_{\Omega_{\mathbf{Q}}}(R(G), J^*(\overline{\mathbf{Q}}))}{Hom_{\Omega_K}(R(G), (\overline{\mathbf{Q}})^*) \cdot Det(\prod_p (\mathbf{Z}_p[G])^*)}.$$

We shall hardly need this Hom-description but, in passing, we remark that a Hom-description for $\pi_*[A, \phi, B]$ is given by the idèlic-valued function which is given at the component above $p$ by any Hom-description representative of $[A \otimes \mathbf{Z}_p, \phi, B \otimes \mathbf{Z}_p] \in K_0(\mathbf{Z}_p[G]; \mathbf{Q}_p)$ and is trivial at Archimedean primes.

**2.4** *The case when $G = Q_8$*

When $G = Q_8$, the quaternion group of order eight, all the complex representations have rational characters so that the isomorphism of §2.3 simplifies to the form

$$K_0(\mathbf{Z}_p[Q_8]; \mathbf{Q}_p) \cong \frac{Hom(R(Q_8), \mathbf{Q}_p^*)}{Det(\mathbf{Z}_p[Q_8]^*)}.$$

Since determinantal functions on $R(Q_8)$ take values in $\mathbf{Z}_p^*$ there is a short exact sequence of the form

$$0 \longrightarrow \frac{Hom(R(Q_8), (\mathbf{Z}_p)^*)}{Det(\mathbf{Z}_p[Q_8]^*)} \longrightarrow \frac{Hom(R(Q_8), (\mathbf{Q}_p)^*)}{Det(\mathbf{Z}_p[Q_8]^*)} \longrightarrow Hom(R(Q_8), \mathbf{Z}) \longrightarrow 0.$$

From this sequence we see that there is an isomorphism of the form

$$TorsK_0(\mathbf{Z}_p[Q_8]; \mathbf{Q}_p) \cong \frac{Hom(R(Q_8), \mathbf{Z}_p^*)}{Det(\mathbf{Z}_p[Q_8]^*)}.$$

When $p$ is an odd prime every function in $Hom(R(Q_8), \mathbf{Z}_p^*)$ is a determinantal function because in this case $\mathbf{Z}_p[Q_8]$ is a maximal order in $\mathbf{Q}_p[Q_8]$ [6]. Hence, if $p$ is odd then $TorsK_0(\mathbf{Z}_p[Q_8]; \mathbf{Q}_p) = 0$.

Now let us examine the case when $p = 2$. Let $1, \chi_1, \chi_2, \chi_1\chi_2$ denote the four one-dimensional representations of $Q_8$ and $Q_8^{ab}$. To be precise, if

$$Q_8 = \{x, y \mid x^2 = y^2, y^4 = 1, xyx = y\}$$

set $\chi_1(x) = -1 = \chi_2(y)$ and $\chi_1(y) = 1 = \chi_2(x)$. The following result is taken from ([16] Proposition 2.5.37 and Corollaries 2.5.38, 2.5.39)

**Proposition 2.5**

(i) There is an isomorphism of the form

$$\lambda : TorsK_0(\mathbf{Z}_2[Q_8^{ab}]; \mathbf{Q}_2) \cong \frac{Hom(R(Q_8^{ab}), \mathbf{Z}_2^*)}{Det(\mathbf{Z}_2[Q_8^{ab}]^*)} \xrightarrow{\cong} (\mathbf{Z}/4)^* = \{\pm 1\}$$

given by

$$\lambda[f] = f(1 + \chi_1 + \chi_2 + \chi_1\chi_2) \quad (modulo\ 4).$$

(ii) The natural maps yield an isomorphism of the form

$$TorsK_0(\mathbf{Z}[Q_8]; \mathbf{Q}) \xrightarrow{\cong} \mathcal{CL}(\mathbf{Z}[Q_8]) \oplus TorsK_0(\mathbf{Z}[Q_8^{ab}]; \mathbf{Q}) \cong \mathbf{Z}/2 \times \mathbf{Z}/2.$$

(iii) $TorsK_0(\mathbf{Z}[\mathbf{Z}/2]; \mathbf{Q}) = 0$.

**Proposition 2.6**

(i) There is an isomorphism of the form

$$K_0(\mathbf{Z}_2[Q_8^{ab}]; \mathbf{Q}_2) \cong Hom(R(Q_8^{ab}), \mathbf{Z}) \times (\mathbf{Z}/4)^*.$$

(ii) Suppose that $H \subset Q_8^{ab}$ has order one or two. Then, under the isomorphism of (i), the canonical homomorphism

$$Ind_H^{Q_8^{ab}} : K_0(\mathbf{Z}_2[H]; \mathbf{Q}_2) \cong \frac{Hom(R(H), \mathbf{Q}_2^*)}{Hom(R(H), \mathbf{Z}_2^*)} \longrightarrow K_0(\mathbf{Z}_2[Q_8^{ab}]; \mathbf{Q}_2)$$

has the form $[f] \mapsto (\chi \mapsto (v_2(f(Res_H^{Q_8^{ab}}(\chi))), 1$ modulo $4)$ where $v_2$ denotes the 2-adic valuation.

**Proof**

For part (i) send the class of a function $f : R(Q_8^{ab}) \longrightarrow \mathbf{Q}_2^*$ to

$$(\chi \mapsto v_2(f(\chi)), f(1 + \chi_1 + \chi_2 + \chi_1\chi_2)2^{-v_2(f(1+\chi_1+\chi_2+\chi_1\chi_2))} \text{ modulo } 4).$$

This an isomorphism by the discussion of §2.4 and Proposition 2.5(i). Since $Ind_H^{Q_8^{ab}}$ is induced in terms of the Hom-description by $Res_H^{Q_8^{ab}}$, part (ii) follows from part (i). □

# 3 Burns invariant $T\Omega^{loc}(N/\mathbf{Q}, 1)$ for quaternion fields

**3.1** Let $N/\mathbf{Q}$ be a quaternion extension of number fields with Galois group $G(N/\mathbf{Q})$ and with biquadratic subfield, $E$. In this case the equivariant Tamagawa number, defined in ([3] §3), is a canonical Galois structure invariant of the form

$$T\Omega^{loc}(N/\mathbf{Q}, 1) \in TorsK_0(\mathbf{Z}[G(N/\mathbf{Q})]; \mathbf{Q}) \cong \mathbf{Z}/2 \times \mathbf{Z}/2.$$

For a general Galois extension of number fields this invariant may only be defineable in the relative K-group, $K_0(\mathbf{Z}[G(N/\mathbf{Q})]; \mathbf{R})$, but for quaternion fields the Stark conjecture holds (cf. [18]) in which case it is shown in [3] that $T\Omega^{loc}(N/\mathbf{Q}, 1)$ is an element of finite order in $K_0(\mathbf{Z}[G(N/\mathbf{Q})]; \mathbf{Q})$.

The following is a particular case of a conjecture of Burn's [3], which is rather remarkable since $T\Omega^{loc}(N/\mathbf{Q}, 1)$ has a definition which looks far from trivial:

**Conjecture 3.2**

$$0 = T\Omega^{loc}(N/\mathbf{Q}, 1) \in TorsK_0(\mathbf{Z}[G(N/\mathbf{Q})]; \mathbf{Q}).$$

The following result will considerably simplify the calculation of $T\Omega^{loc}(N/\mathbf{Q}, 1)$.

**Lemma 3.3**

The image of $T\Omega^{loc}(N/\mathbf{Q}, 1)$ under the isomorphism of Proposition 2.5(ii) is equal to

$$(0, T\Omega^{loc}(E/\mathbf{Q}, 1)) \in \mathcal{CL}(\mathbf{Z}[Q_8]) \oplus TorsK_0(\mathbf{Z}[Q_8^{ab}]; \mathbf{Q}) \cong \mathbf{Z}/2 \times \mathbf{Z}/2.$$

**Proof**

It is shown in in [3] that the image of $T\Omega^{loc}(N/\mathbf{Q}, 1)$ in the class group is equal to $\Omega(N/\mathbf{Q}, 2) - W_{N/\mathbf{Q}}$, the difference of the second Chinburg invariant and the Cassou-Noguès-Fröhlich class. However this difference is shown to be trivial for all quaternion fields in [7]. This shows that the first component of the image vanishes. By naturality of $T\Omega^{loc}(N/\mathbf{Q}, 1)$ with respect to passage to Galois subextensions (cf. [1]) the second component is equal to $T\Omega^{loc}(E/\mathbf{Q}, 1)$. □

**Definition 3.4** *The invariant $T\Omega^{loc}(E/\mathbf{Q}, 1)$*

Let $E = \mathbf{Q}(\sqrt{d_1}, \sqrt{d_2})$ be the biquadratic subfield of the quaternion field $N/\mathbf{Q}$, as in §3.1. By virtue of Lemma 3.3, we wish to calculate the invariant $T\Omega^{loc}(E/\mathbf{Q}, 1)$ of ([3] Proposition 3.3 p.15). In general this invariant is defined as a difference

$$T\Omega^{loc}(E/\mathbf{Q},1) = T\Omega^{loc}(E/\mathbf{Q},1,\lambda) - T(\lambda) \in K_0(\mathbf{Z}[G(E/\mathbf{Q})],\mathbf{Q})$$

depending on an element, $\lambda \in \mathbf{R}[G(E/\mathbf{Q})]^*$. However, when all units in the real group-ring are reduced norms we may choose $\lambda = 1$ and then $T(1)$ is trivial.

The definition of $T\Omega^{loc}(E/\mathbf{Q},1,\lambda)$ requires the choice of a finite, Galois invariant set of places of $E$, denoted by $S$ and containing all the infinite places $S_\infty(E)$ together with all finite places which ramify in $E/\mathbf{Q}$. Then, in the notation of ([3] p.13),

$$T\Omega^{loc}(E/\mathbf{Q},1,\lambda) = T_S(\mathcal{L},\lambda) - \sum_{p\in S_f} Ind_{G(E_w/\mathbf{Q}_p)}^{G(E/\mathbf{Q})}(\chi_w^\bullet(\mathcal{L},V_w)) + \sum_{p\in S_f} \hat{\delta}_1(\epsilon'_p(0)).$$

Here $w$ is a chosen prime of $E$ above $p$ and $\mathcal{L}$ is a lattice of which more presently. Each of the sums is taken over $S_f$, which is the set of rational primes $p$ lying below the finite primes of $E$ in $S$. The terms $\hat{\delta}_1(\epsilon'_p(0))$ of [3] are elements lying in the image of the map $\delta_1$ of §2.1.

**3.5** If $\chi : G(E/\mathbf{Q}) \longrightarrow \{\pm 1\}$ is a one-dimensional representation a Hom-representative of $\hat{\delta}_1(\epsilon'_p(0)) \in K_0(\mathbf{Z}[G(E/\mathbf{Q})],\mathbf{Q})$ has the same component at each prime $l$ given by the function

$$\chi \mapsto (\frac{|G(E_w/\mathbf{Q}_p)|}{|I(E_w/\mathbf{Q}_p)|})^{-dim(\chi^{G(E_w/\mathbf{Q}_p)})} det(1 - Frob_w^{-1}|(\chi^{I(E_w/\mathbf{Q}_p)}/\chi^{G(E_w/\mathbf{Q}_p)}))$$

where $I(E_w/\mathbf{Q}_p) \subseteq G(E_w/\mathbf{Q}_p)$ is the inertia subgroup at $p$. Here $\chi^H$ denotes the subrepresentation given by the $H$-fixed points for $H \subseteq Q_8^{ab}$. Hence the Hom-representative of $\hat{\delta}_1(\epsilon'_p(0))$ has the form $\chi \mapsto 2^{\alpha_p(\chi)}$. Since $2 \in \mathbf{Z}_l^*$ for all odd primes, a Hom-representative of $\sum_{p\in S_f} \hat{\delta}_1(\epsilon'_p(0))$ is given by the function which is given by

$$[\chi \mapsto 2^{\sum_{p\in S_f} \alpha_p(\chi)}] \in \frac{Hom(R(Q_8^{ab}),\mathbf{Q}_2^*)}{Det(\mathbf{Z}_2[Q_8^{ab}]^*)}$$

in the 2-adic coordinate and is trivial at all odd-primary components.

Combining this discussion with Proposition 2.6(i) yields the following result.

**Corollary 3.6**

*The torsion component of the image of $\sum_{p\in S_f} \hat{\delta}_1(\epsilon'_p(0))$ in*

$$K_0(\mathbf{Z}_2[Q_8^{ab}];\mathbf{Q}_2) \cong Hom(R(Q_8^{ab}),\mathbf{Z}) \times (\mathbf{Z}/4)^*$$

*is trivial.*

**3.7** *The element $T_S(\mathcal{L}, 1)$*

The term, $T_S(\mathcal{L}, \lambda)$ (with $\lambda = 1$ in our case) is constructed ([3] p.10) in terms of a free $\mathbf{Z}[G(E/\mathbf{Q})]$-sublattice $\mathcal{L}$ of the integers, $\mathcal{O}_E$, of $E$ so that $\mathcal{L} = \mathbf{Z}[G(E/\mathbf{Q})] < \gamma >$ for some $\gamma \in \mathcal{O}_E$. In general, $\mathcal{L}$ would merely be a projective module but finitely generated, projective $\mathbf{Z}[G(E/\mathbf{Q})]$-modules are all free.

Let $\mathbf{Z} < 2\pi i >$ denote a copy of the integers on which complex conjugation acts by minus the identity. Form the chain complex of perfect $\mathbf{Z}[G(E/\mathbf{Q})]$-modules

$$H_B^\bullet(\mathcal{L}): \ 0 \longrightarrow \mathbf{Z} < 2\pi i > \otimes \mathbf{Z}[G(E/\mathbf{Q})] \xrightarrow{(2,0)} \mathbf{Z} < 2\pi i > \otimes \mathbf{Z}[G(E/\mathbf{Q})] \oplus \mathcal{L}$$

$$\xrightarrow{0} \mathbf{Z} < 2\pi i > \otimes \mathbf{Z}[G(E/\mathbf{Q})] \longrightarrow 0$$

in which the non-zero modules are, from left to right, in degrees $-1, 0, 1$, respectively. Hence $H^0 = \mathcal{L}$ and $H^1 = \mathbf{Z} < 2\pi i > \otimes \mathbf{Z}[G(E/\mathbf{Q})]$. The class, $T_S(\mathcal{L}, 1) \in K_0(\mathbf{Z}[G(E/\mathbf{Q})], \mathbf{Q})]$, is the element associated to this perfect complex together with a rational homology isomorphism defined in the manner of §2.2. There is a rational isomorphism

$$\Theta : \mathcal{L} \otimes \mathbf{Q} \xrightarrow{\cong} \mathbf{Z} < 2\pi i > \otimes \mathbf{Q}[G(E/\mathbf{Q})]$$

which, for $\gamma \in \mathcal{L}$ such that $\gamma \otimes 1$ is a generator for the module $\mathcal{L} \otimes \mathbf{Q}$, is given by the formula

$$\Theta(\gamma \otimes 1) = (2\pi i) \otimes \sum_{g \in G(E/\mathbf{Q})], \mathbf{Q})} g(\gamma) g^{-1}.$$

The rational isomorphism used in ([3] p.6), obtained by using $\Theta^{-1} : H^{od} \otimes \mathbf{Q} \xrightarrow{\cong} H^{ev} \otimes \mathbf{Q}$ in §2.2, is equal to

$$\oplus_{i=1,2} \mathbf{Z} < 2\pi i > \otimes \mathbf{Q}[G(E/\mathbf{Q})] \xrightarrow{\cong} \mathbf{Z} < 2\pi i > \otimes \mathbf{Q}[G(E/\mathbf{Q})] \oplus \mathcal{L} \otimes \mathbf{Q}$$

given by

$$\mathcal{E}_S^{-1} \begin{pmatrix} 2 & 0 \\ 0 & \Theta^{-1} \end{pmatrix}$$

Here $\mathcal{E}_S \in \mathbf{R}[G(E/\mathbf{Q})]^*$ is defined in ([3] p.6) and will be recalled in detail below. Therefore a Hom-representative for $T_S(\mathcal{L}, 1)$ is given by the function whose coordinate at each prime $p$ is given on a one-dimensional representation, $\chi$, by

$$\chi \mapsto 2\chi(\mathcal{E}_S^{-1}) det(\sum_{g \in G(E/\mathbf{Q})} g(\gamma) \chi(g^{-1}))^{-1} = 2\chi(\mathcal{E}_S^{-1})(\gamma \mid \chi)^{-1},$$

where $(\gamma \mid \chi)$ is the resolvent (see [5]) associated to $\gamma \in E$. Finally, for $\chi : G(E/\mathbf{Q}) \longrightarrow \{\pm 1\}$, $\chi(\mathcal{E}_S^{-1})$ is given by the formula (see [3] §1.2.2.1(proof))

$$\chi(\mathcal{E}_S) = \frac{L^*(1, \chi)}{L^*(0, \chi)} \prod_{p \in S_f} det(1 - p^{-1} Frob_w^{-1} \mid \chi^{I(E_w/\mathbf{Q}_p)})$$

where $L^*(s_0, \chi)$ denotes the leading coefficient in the Laurent series for the Artin L-function, $L(s,\chi)$, about $s = s_0$.

**Lemma 3.8**

*In §3.7*
$$\frac{L^*(1,\chi)}{L^*(0,\chi)} = 2f(\chi)^{-1/2} = 2f(\chi)^{-1/2}$$
*for each one-dimensional representation, $\chi : G(E/\mathbf{Q}) \longrightarrow \{\pm 1\}$.*

**Proof**

In order to compute the quotient of the leading coefficients we recall that there is a functional equation ([10], [13] p.253) of the form

$$f(\chi)^{s/2}\pi^{-s/2}\Gamma(s/2)L(s,\chi)$$
$$= f(\chi)^{(1-s)/2}\pi^{-(1-s)/2}\Gamma((1-s)/2)L(1-s,\chi)$$

in which $f(\chi)$ denotes the Artin conductor of $\chi$. In the neighbourhood of $s = 1$ recall that when $\chi = 1$ the L-function has a simple pole (as is seen from the Analytic Class Number Formula (see e.g. [13] p.248)) while $\Gamma(s/2)$ is continuous and non-zero in the vicinity of $s = 1$. At $s = 0$ the behaviour is the other way round. When $\chi$ is non-trivial the L-function tends to a non-zero limit at $s = 1$. Therefore re-writing the function equation as

$$f(\chi)^{s/2}\pi^{-s/2}\frac{L(s,\chi)}{\Gamma((1-s)/2)} = f(\chi)^{(1-s)/2}\pi^{-(1-s)/2}\frac{L(1-s,\chi)}{\Gamma(s/2)}.$$

Consider the case when $\chi = 1$. Hence

$$\lim_{t \to 1}(t-1)L(t,1) = \lim_{s \to 0}(1-s-1)L(1-s,1) = -\lim_{s \to 0} sL(1-s,1)$$

exists so that near $s = 0$

$$L(1-s,1) = L^*(1,1)/s + \sum_{j=0}^{\infty} z_j s^j.$$

Similarly it is well known ([9] p.8 and p.491) that near $s = 0$

$$\Gamma(s/2) = 2/s + \sum_{j=0}^{\infty} y_j s^j.$$

Letting $s$ tend to zero in the functional equation yields

$$\pi^{-1/2}L^*(0,1) = \frac{L^*(0,1)}{\Gamma(1/2)} = \pi^{-1/2}\frac{L^*(1,1)}{2}.$$

Hence
$$\frac{L^*(1,1)}{L^*(0,1)} = 2.$$

When $\chi$ is non-trivial the limit of $L(s,\chi)$ as $s$ tends to one exists and is non-zero, since $L(s,\chi)$ is the quotient of two zeta functions, and so
$$L^*(1,\chi) = \lim_{s \to 1} L(s,\chi) = \lim_{s \to 0} L(1-s,\chi).$$

Therefore near $s = 0$
$$L(s,\chi) = L^*(0,\chi)s + \sum_{j=2}^{\infty} u_j s^j$$

and dividing the functional equation by $s$ and letting $s$ tend to zero yields
$$\pi^{-1/2} L^*(0,\chi) = f(\chi)^{1/2} \pi^{-1/2} \frac{L^*(1,\chi)}{2}$$

so that
$$\frac{L^*(1,\chi)}{L^*(0,\chi)} = 2f(\chi)^{-1/2},$$

as required. $\square$

**Corollary 3.9**

Hence, for each one-dimensional representation $\chi : G(E/\mathbf{Q}) \longrightarrow \{\pm 1\}$, a Hom-representative for $T_S(\mathcal{L}, 1)$ of §3.7 is given by the function whose coordinate at each prime is

$$\chi \mapsto \frac{f(\chi)^{1/2}}{(\gamma \mid \chi) \prod_{p \in S_f} det(1 - p^{-1} Frob_w^{-1} \mid \chi^{I(E_w/\mathbf{Q}_p)})}.$$

**Proposition 3.10**

By Proposition 2.6(ii) the image of $Ind_{G(E_w/\mathbf{Q}_p)}^{G(E/\mathbf{Q})}(\chi_w^\bullet(\mathcal{L}), V_w)$ in

$$K_0(\mathbf{Z}_2[Q_8^{ab}]; \mathbf{Q}_2) \cong Hom(R(Q_8^{ab}), \mathbf{Z}) \times (\mathbf{Z}/4)^*$$

has trivial torsion component unless $G(E_w/\mathbf{Q}_p) = G(E/\mathbf{Q})$.

# 4 The local terms $Ind_{G(E_w/\mathbf{Q}_p)}^{G(E/\mathbf{Q})}(\chi_w^\bullet(\mathcal{L}, V_w))$

**Definition 4.1** The terms in the middle expression for $T\Omega^{loc}(E/\mathbf{Q}, 1, \lambda)$ in Definition 3.4 are the images of elements, $\chi_w^\bullet(\mathcal{L}, V_w) \in Tors K_0(\mathbf{Z}[G(E_w/\mathbf{Q}_p)], \mathbf{Q})$ under the canonical homomorphism

$$Ind_{G(E_w/\mathbf{Q}_p)}^{G(E/\mathbf{Q})} : K_0(\mathbf{Z}[G(E_w/\mathbf{Q}_p)], \mathbf{Q}) \longrightarrow Tors K_0(\mathbf{Z}[G(E/\mathbf{Q})], \mathbf{Q})$$

where $G(E_w/\mathbf{Q}_p)$ denotes the decomposition group at $p$. Hence, by Proposition 3.10, the torsion part of the local contribution, $Ind_{G(E_w/\mathbf{Q}_p)}^{G(E/\mathbf{Q})}(\chi_w^\bullet(\mathcal{L}, V_w))$, vanishes unless $G(E_w/\mathbf{Q}_p) = G(E/\mathbf{Q})$. In this section we shall evaluate these local torsion contributions when $p$ is odd and $G(E_w/\mathbf{Q}_p) = G(E/\mathbf{Q})$. In this case $p$ is ramified, so that $p \in S_f$, and the inertia group $I(E_w/\mathbf{Q}_p)$ has order two. In this case, let us describe the class $\chi_w^\bullet(\mathcal{L}, V_w)$. Choose $w \in S$ over $p$. Let $\mathcal{L}_w \subset E_w$ denote the $w$-adic completion of $\mathcal{L}$. According to [3] we may choose any free $\mathcal{L}$ such that the $p$-adic exponential is defined and gives an isomorphism,

$$exp: \mathcal{L}_w \xrightarrow{\cong} 1 + \mathcal{L}_w.$$

Since an odd ramified prime is tamely ramified, $\mathcal{O}_{E_w}$ is a free $\mathbf{Z}_p[G(E_w/\mathbf{Q}_p)]$-module and we may choose an integer $m = m(p) \geq 1$ and arrange that $p^m \mathcal{O}_{E_w} = \mathcal{L}_w$ (cf. [6], [14]).

Now suppose that the 2-extension

$$E_w^* \longrightarrow A \longrightarrow B \longrightarrow \mathbf{Z}$$

respresents the canonical local fundamental class of class field theory in

$$Ext^2_{\mathbf{Z}[G(E_w/\mathbf{Q}_p)]}(\mathbf{Z}, E_w^*) = H^2(G(E_w/\mathbf{Q}_p); E_w^*)$$

([12] p.168). Form the 2-extension

$$E_w^*/(1 + p^m \mathcal{O}_{E_w}) \longrightarrow A/(1 + p^m \mathcal{O}_{E_w}) \longrightarrow B \longrightarrow \mathbf{Z}.$$

Choose a complex

$$P^*: \; 0 \longrightarrow P^{-2} \xrightarrow{d_{-2}} P^{-1} \xrightarrow{d_{-1}} P^0 \longrightarrow 0$$

of finitely generated, projective $\mathbf{Z}[G(E_w/\mathbf{Q}_p)]$-modules which is quasi-isomorphic to

$$0 \longrightarrow A/(1 + p^m \mathcal{O}_{E_w}) \longrightarrow B \longrightarrow 0$$

in which $B$ is in degree zero. Hence there are canonical isomorphisms of the form

$$H^{-1}(P^*) \cong E_w^*/(1 + p^m \mathcal{O}_{E_w}), \quad H^0(P^*) \cong \mathbf{Z}.$$

Therefore the valuation gives a rational isomorphism of the form

$$val: H^{-1}(P^*) \otimes \mathbf{Q} \xrightarrow{\cong} H^0(P^*) \otimes \mathbf{Q}$$

and by the construction of §2.2 we obtain

$$[P^{-1}, val, P^{-2} \oplus P^0] = \chi_w^\bullet(\mathcal{L}, V_w) \in K_0(\mathbf{Z}[G(E_w/\mathbf{Q}_p)], \mathbf{Q}).$$

**4.2** *Computing the local terms at odd primes*

We continue with the situation of Definition 4.1. In $G(E/\mathbf{Q}) = G(E_w/\mathbf{Q}_p)$ let $a_p$ denote the non-trivial element of the inertia group at $p$ and let $b_p$ be an element acting like the Frobenius of the residue extension, $\mathbf{F}_{p^2}/\mathbf{F}_p$.

In order to compute the local term we are first going study the related chain complex

$$0 \longrightarrow A/U^1_{E_w} \longrightarrow B \longrightarrow 0$$

in which we quotient by the principal units, $U^1_{E_w}$, rather than by $1+p^m \mathcal{L}_w$. Since $p$ is tamely ramified, $U^1_{E_w}$ is a cohomologically trivial $\mathbf{Z}[G(E_w/\mathbf{Q}_p)]$-module.

If, as in Definition 4.1,

$$E^*_w \longrightarrow A \xrightarrow{d} B \longrightarrow \mathbf{Z}$$

respresents the canonical local fundamental class of class field theory then minus the canical class is given by pulling this back along $(-1) : \mathbf{Z} \longrightarrow \mathbf{Z}$. Hence minus the fundamental class is represented by

$$E^*_w \longrightarrow A \xrightarrow{-d} B \longrightarrow \mathbf{Z}.$$

Since $p$ is tamely ramified we have a complex of free $\mathbf{Z}[G(E_w/\mathbf{Q}_p)]$-modules which is quasi-isomorphic to minus the local fundamental class divided by $U^1_{E_w}$

$$E^*_w/(U^1_{E_w}) \longrightarrow A/(U^1_{E_w}) \xrightarrow{-d} B \longrightarrow \mathbf{Z}$$

which was first constructed in [4] and is constructed in a slightly different manner in [15] pp.319-329). This quasi-isomorphic complex of free modules has the form

$$0 \longrightarrow \mathbf{Z}[G(E/\mathbf{Q})] <w> \xrightarrow{\lambda} \mathbf{Z}[G(E/\mathbf{Q})] <z_1> \oplus \mathbf{Z}[G(E/\mathbf{Q})] <z_2>$$

$$\xrightarrow{\phi} \mathbf{Z}[G(E/\mathbf{Q})] <t> \longrightarrow 0$$

where

$$\lambda(w) = (b_p((p+1)/2 + ((p-1)/2)a_p) - 1)z_1 - (a_p - 1)z_2$$

$$\phi(z_1) = (a_p - 1)t, \ \phi(z_2) = (b_p - 1)t.$$

If the left-hand module is in dimension minus two then the only non-trivial cohomology of this complex is given by

$$H^{-1} \cong E^*_w/U^1_{E_w} \cong \mathbf{F}^*_{p^2} \times \mathbf{Z}, \quad H^0 \cong \mathbf{Z}$$

where the first isomorphism sends $T = (1+b_p)z_2 - (1+a_p)z_1$ to $(\xi, \pi)$, by ([15] Lemma 7.1.55) – $\pi$ being a prime of $E_w$ – and the second isomorphism is induced by the augmentation homomorphism. The valuation on $E_w$, normalised to ensure that $val(\pi) = 1$, induces a canonical isomorphism, $H^{-1} \otimes \mathbf{Q} \longrightarrow H^0 \otimes \mathbf{Q}$.

Therefore the complex

$$0 \longrightarrow \mathbf{Z}[G(E/\mathbf{Q})] <w> \xrightarrow{\lambda} \mathbf{Z}[G(E/\mathbf{Q})] <z_1> \oplus \mathbf{Z}[G(E/\mathbf{Q})] <z_2>$$

$$\xrightarrow{-\phi} \mathbf{Z}[G(E/\mathbf{Q})] <t> \longrightarrow 0$$

is quasi-isomorphic to the quotient of the local fundamental class by $U^1_{E_w}$. Let us apply the construction of §2.2 to this complex and rational cohomology isomorphism. We must choose splittings

$$0 \longrightarrow im(-\phi) \otimes \mathbf{Q} \xrightarrow{i_0} \mathbf{Q}[G(E/\mathbf{Q})] <t> \xleftarrow{\rho_0} H^0 \otimes \mathbf{Q} \longrightarrow 0,$$

$$0 \longrightarrow ker(-\phi) \otimes \mathbf{Q} \xrightarrow{i_1} \mathbf{Q}[G(E/\mathbf{Q})] <z_1, z_2> \xleftarrow{\rho_1} im(-\phi) \otimes \mathbf{Q} \longrightarrow 0,$$

$$0 \longrightarrow \mathbf{Q}[G(E/\mathbf{Q})] <w> \xrightarrow{\lambda} ker(-\phi) \otimes \mathbf{Q} \xleftarrow{\rho_2} H^{-1} \otimes \mathbf{Q} \longrightarrow 0.$$

Here $i_j$ is induced by an inclusion map.

The recipe of §2.2 then makes use of the following string of isomorphisms

$$\mathbf{Q}[G(E/\mathbf{Q})] <z_1, z_2>$$

$$\xleftarrow{i_1+\rho_1} ker(-\phi) \otimes \mathbf{Q} \oplus im(-\phi) \otimes \mathbf{Q}$$

$$\xleftarrow{(i_2+\rho_2, 1)} im(\lambda) \otimes \mathbf{Q} \oplus im(-\phi) \otimes \mathbf{Q} \oplus H^{-1} \otimes \mathbf{Q}$$

$$\xleftarrow{(\lambda, 1, (val)^{-1})} \mathbf{Q}[G(E/\mathbf{Q})] <w> \oplus im(-\phi) \otimes \mathbf{Q} \oplus H^0 \otimes \mathbf{Q}$$

$$\xrightarrow{(1, i_0+\rho_0)} \mathbf{Q}[G(E/\mathbf{Q})] <w, t>.$$

The decomposition of the group-ring $\mathbf{Q}[G(E/\mathbf{Q})]$ into the product of one copy of $\mathbf{Q}$ for each one-dimensional representation is induced by evaluating at the corresponding representation. Hence we can produce the function, $f \in Hom(R(G(E/\mathbf{Q}), \mathbf{Q}^*)$, which gives the $p$-component (for every $p$) of the Hom-representative for the element of $K_0(\mathbf{Z}[G(E/\mathbf{Q})], \mathbf{Q})$ associated with the above isomorphism data idempotent by idempotent. In other words, if $\sigma_\chi \in \mathbf{Q}[G(E/\mathbf{Q})]$ is the idempotent corresponding to $\chi : G(E/\mathbf{Q}) \longrightarrow \{\pm 1\}$ then $f(\chi)$ is given by the determinant of the matrix representing the restriction of the above isomoprhism

$$\mathbf{Q} <\sigma_\chi z_1, \sigma_\chi z_2> \longrightarrow \mathbf{Q} <\sigma_\chi w, \sigma_\chi t>.$$

Furthermore, the splittings, $\rho_i$, are determined by the splittings, $\sigma_\chi \rho_i$.

We begin with $\chi = 1$, the one-dimensional trivial representation, so that $\sigma_1 = (a_p + 1)(b_p + 1)/4$. In this case $\sigma_1 im(\phi) \otimes \mathbf{Q} = 0$ and $\rho_0(\sigma_1 t) = 1 \otimes 1 = val(\pi \otimes 1) \in H_0 \otimes \mathbf{Q}$. Since $\sigma_1 \lambda(w) = (p-1)\sigma_1 z_1$ and we may choose $\rho_2(\sigma_1(\pi \otimes 1)) = \sigma_1 T = 2\sigma_1(z_2 - z_1)$ we obtain

$$f(1) = det \begin{pmatrix} 1/(p-1) & 1/(p-1) \\ 0 & 1/2 \end{pmatrix} = 1/(2p-2) \in \mathbf{Q}^*.$$

When $\chi$ is a non-trivial one-dimensional representation then $\sigma_\chi = (1 + \chi(a_p)a_p)(1 + \chi(b_p)b_p)/4$ and $\sigma_\chi H^{-1} \otimes \mathbf{Q} = 0 = \sigma_\chi H^0 \otimes \mathbf{Q}$.

Suppose that $\chi(a_p) = -1 = \chi(b_p)$. Then $(-\phi)(\sigma_\chi z_1) = 2\sigma_\chi t = (-\phi)(\sigma_\chi z_2)$ so that we may choose $\rho_1(\sigma_\chi t) = \sigma_\chi z_1/2$ since $(a_p - 1)(-\sigma_\chi z_1/2) = \sigma_\chi t =$. Also $\lambda(\sigma_\chi w) = 2\sigma_\chi (z_2 - z_1)$. Therefore, in this case, we obtain

$$f(\chi) = det \begin{pmatrix} 0 & 1/2 \\ 2 & 2 \end{pmatrix} = -1 \in \mathbf{Q}^*.$$

If $\chi(a_p) = 1$, $\chi(b_p) = -1$ then $\lambda(\sigma_\chi w) = -(p+1)\sigma_\chi z_1$ and we may choose $\rho_1(\sigma_\chi t) = \sigma_\chi z_2/2$. Therefore we obtain

$$f(\chi) = det \begin{pmatrix} -1/(p+1) & 0 \\ 0 & 2 \end{pmatrix} = -2/(p+1) \in \mathbf{Q}^*.$$

Finally, if $\chi(a_p) = -1$, $\chi(b_p) = 1$ then $\lambda(\sigma_\chi w) = 2\sigma_\chi z_2$ and we may choose $\rho_1(\sigma_\chi t) = \sigma_\chi z_1/2$. Therefore

$$f(\chi) = det \begin{pmatrix} 0 & 1/2 \\ 2 & 0 \end{pmatrix} = -1 \in \mathbf{Q}^*.$$

One may easily verify that this function may be written as the $\mathbf{Q}^*$-valued function sending a one-dimensional $\chi$ to

$$(-1)^{codim(\chi^{G(E_w/\mathbf{Q}_p)})} \frac{(\frac{|G(E_w/\mathbf{Q}_p)|}{|I(E_w/\mathbf{Q}_p)|})^{-dim(\chi^{G(E_w/\mathbf{Q}_p)})} det(1 - Frob_w^{-1}|(\chi^{I(E_w/\mathbf{Q}_p)}/\chi^{G(E_w/\mathbf{Q}_p)}))}{p^{dim(\chi^{I(E_w/\mathbf{Q}_p)})} det(1 - p^{-1} Frob_w^{-1} \mid \chi^{I(E_w/\mathbf{Q}_p)})}.$$

We still have to evaluate these local terms using the correct lattice.

The categorical derivation of the localisation sequence in [11] shows that $K_0(\mathbf{Z}[G(E_w/\mathbf{Q}_p)], \mathbf{Q}) \cong K_0T(\mathbf{Z}[G(E_w/\mathbf{Q}_p)])$, the Grothendieck group of the category of finite $\mathbf{Z}[G(E_w/\mathbf{Q}_p)]$-modules of finite projective dimension. The module $U^1_{E_w}/(1 + p^m \mathcal{O}_{E_w})$ lies in this category. To see this observe that each quotient in the level filtration $(1 + p^{n-1}\mathcal{O}_{E_w})/(1 + p^n\mathcal{O}_{E_w})$ is isomorphic to the residue field $\mathbf{F}_{p^2}$, which has a free resolution of the form

$$0 \longrightarrow \mathbf{Z}[G(\mathcal{O}_{E_w}/\mathbf{Q}_p)] \stackrel{((p+1)/2+(1/2)(p-1)a_p)\cdot -}{\longrightarrow} \mathbf{Z}[G(\mathcal{O}_{E_w}/\mathbf{Q}_p)] \longrightarrow \mathbf{F}_{p^2} \longrightarrow 0.$$

Here the second map sends 1 to a normal basis element of the residue field. To verify exactness we note that $(p+1)/2 + (1/2)(p-1)a_p$ acts on the normal basis element as multiplcation by $p$ so that the composition of the two maps is zero. Since $(p+1)/2 + (1/2)(p-1)a_p$ maps to 1 or $p$ under the one-dimensional representations of $G(\mathcal{O}_{E_w}/\mathbf{Q}_p)$ one sees that the left-hand map is injective. To complete the verification of exactness we must show that the ideal, $< p, (a_p-1) > \triangleleft \mathbf{Z}[G(\mathcal{O}_{E_w}/\mathbf{Q}_p)]$ equals the principal ideal generated by $(p+1)/2 + (1/2)(p-1)a_p$. It is clear that $(p+1)/2 + (1/2)(p-1)a_p \in < p, (a_p - 1) >$ and to see the reverse inclusion it suffices to observe that

$$((p+1)/2 + (1/2)(p-1)a_p) - a_p((p+1)/2 + (1/2)(p-1)a_p)$$

$$= (p+1)/2 + (1/2)(p-1)a_p - (1/2)(p+1)a_p - (1/2)(p-1)$$

$$= 1 - a_p.$$

This shows that $[\mathbf{F}_{p^2}] \in K_0T(\mathbf{Z}[G(E_w/\mathbf{Q}_p)]) \cong K_0(\mathbf{Z}[G(E_w/\mathbf{Q}_p)], \mathbf{Q})$ has a Hom-description representative given by the function whose coordinate at each prime sends a one-dimensional representation $\chi$ of $G(E_w/\mathbf{Q}_p)$ to

$$\chi((p+1)/2 + (1/2)(p-1)a_p) = \begin{cases} p & \text{if } \chi(a_p) = 1, \\ 1 & \text{otherwise.} \end{cases}$$

if we consider $\mathbf{F}_{p^2}$ to be an *odd-dimensional* homology group. Since $p = p^{-1} \in (\mathbf{Z}/4)^*$, where our final computation takes place, the homological grading dimension of $\mathbf{F}_{p^2}$ will be immaterial.

This discussion shows that the Hom-description representative of $\chi_w^\bullet(\mathcal{L}, V_w)$ differs from the one obtained above, using $U_{E_w}^1$, by $(\chi \mapsto p^{mdim(\chi^{I(E_w/\mathbf{Q}_p)})})^{\pm 1}$. As remarked earlier, for our purposes the sign of the exponent will be immaterial.

**Proposition 4.3**

*Let $p$ be an odd prime such that $G(E/\mathbf{Q}) = G(E_w/\mathbf{Q}_p)$, as in Definition 4.1 and §4.2. Then a Hom-description representative for the local element*

$$\chi_w^\bullet(\mathcal{L}, V_w) \in TorsK_0(\mathbf{Z}[G(E_w/\mathbf{Q}_p)], \mathbf{Q})$$

*is given by the function whose coordinate at each prime on one-dimension representations $\chi$ equals*

$$\epsilon(\chi) \frac{(\frac{|G(E/\mathbf{Q})|}{|I(E_w/\mathbf{Q}_p)|})^{-dim(\chi^{G(E_w/\mathbf{Q}_p)})} det(1 - Frob_w^{-1}|(\chi^{I(E_w/\mathbf{Q}_p)}/\chi^{G(E_w/\mathbf{Q}_p)}))}{p^{1 \pm mdim(\chi^{I(E_w/\mathbf{Q}_p)})} det(1 - p^{-1}Frob_w^{-1} \mid \chi^{I(E_w/\mathbf{Q}_p)})}$$

*where*
$$\epsilon(\chi) = (-1)^{dim(\chi^{I(E_w/\mathbf{Q}_p)}/\chi^{G(E_w/\mathbf{Q}_p)})}.$$

## 5 Conjecture 2.2 for some quaternion fields

**5.1** In this section we shall verify Conjecture 2.2 in the case when the decomposition group at $p = 2$ $|G(N_v/\mathbf{Q}_2)|$ strictly smaller than $Q_8$. The remaining cases correspond to quaternion extensions $N/\mathbf{Q}$ in Cases A,B or C in [7]. These remaining cases seem to be much more difficult and will hopefully be treated in [8].

**Theorem 5.2**
 *Let $N/\mathbf{Q}$ be a quaternion extension of number fields with Galois group $G(N/\mathbf{Q})$. If the decomposition group at $p = 2$ is strictly smaller than $Q_8$ then*
$$T\Omega^{loc}(N/\mathbf{Q}, 1) \in TorsK_0(\mathbf{Z}[G(N/\mathbf{Q})]; \mathbf{Q}) \cong \mathbf{Z}/4)^*$$
*is trivial, as predicted in Conjecture 2.2.*

**Proof**
Let $E$ denote the biquadratic subfield of $N$. By Lemma 3.3 we must show that
$$T\Omega^{loc}(E/\mathbf{Q}, 1)) \in TorsK_0(\mathbf{Z}[G(E/\mathbf{Q})]; \mathbf{Q}) \cong (\mathbf{Z}/4)^*$$
is trivial. By Definition 3.4
$$T\Omega^{loc}(E/\mathbf{Q}, 1) = T_S(\mathcal{L}, \lambda) - \sum_{p \in S_f} Ind_{G(E_w/\mathbf{Q}_p)}^{G(E/\mathbf{Q})}(\chi_w^\bullet(\mathcal{L}, V_w)) + \sum_{p \in S_f} \hat{\delta}_1(\epsilon'_p(0)).$$

Let $\chi_1, \chi_2 : G(E/\mathbf{Q}) \longrightarrow \{\pm 1\}$ be two distinct, non-trivial characters. By Proposition 2.6 and the discussion of §2.4, if $f : R(G(E/\mathbf{Q})) \longrightarrow \mathbf{Q}_2^*$ is the 2-component of a Hom-description representative for $T\Omega^{loc}(E/\mathbf{Q}, 1) \in (\mathbf{Z}/4)^*$, then
$$T\Omega^{loc}(E/\mathbf{Q}, 1) = f(1 + \chi_1 + \chi_2 + \chi_1\chi_2)2^{-v_2(f(1+\chi_1+\chi_2+\chi_1\chi_2))} \text{ (modulo 4)}.$$

By Corollary 3.6 and Proposition 3.10 this remains true if $f$ is replaced by $h$, a Hom-description representative for
$$T_S(\mathcal{L}, 1) - \sum_{p \in S_f, \; G(E_w/\mathbf{Q}_p) = G(E/\mathbf{Q})} \chi_w^\bullet(\mathcal{L}, V_w).$$

Note that, by hypothesis, the primes in the some are all odd.
Let $h_p$ denote the Hom-description respresentative given in Proposition 4.3 for the local term corresponding to an odd prime $p$ in the sum. The factor $\epsilon(1+\chi_1+\chi_2+\chi_1\chi_2) = \epsilon(1)\epsilon(\chi_1)\epsilon(\chi_2)\epsilon(\chi_1\chi_2) = 1$ because there are precisely two distinct, non-trivial characters which are trivial on the inertia group. For the same reason

$$\prod_{\chi=1,\chi_1,\chi_2,\chi_1\chi_2} p^{1\pm mdim(\chi^{I(E_w/\mathbf{Q}_p)})}$$

is a square in $(\mathbf{Z}/4)^*$ and hence trivial. The factors in the numerator are all positive powers of two, which cancel out in $h_p(1+\chi_1+\chi_2+\chi_1\chi_2)2^{-v_2(h_p(1+\chi_1+\chi_2+\chi_1\chi_2))}$. Therefore the torsion component of the two of local terms is equal to

$$\prod_{\substack{p\in S_f,\ G(E_w/\mathbf{Q}_p)=G(E/\mathbf{Q}) \\ \chi=1,\chi_1,\chi_2,\chi_1\chi_2}} det(1-p^{-1}Frob_w^{-1} \mid \chi^{I(E_w/\mathbf{Q}_p)})2^{-v_2(det(1-p^{-1}Frob_w^{-1} \mid \chi^{I(E_w/\mathbf{Q}_p)}))}$$

modulo 4.

By Corollary 3.9 a Hom-description representative for $T_S(\mathcal{L},1)$ of §3.7 is given by the function whose 2-adic coordinate is

$$\chi \mapsto \frac{f(\chi)^{1/2}}{(\gamma \mid \chi) \prod_{p\in S_f} det(1-p^{-1}Frob_w^{-1} \mid \chi^{I(E_w/\mathbf{Q}_p)})}.$$

By [19], the function $(\chi \mapsto \frac{f(\chi)^{1/2}}{(\gamma \mid \chi)}) \in \mathbf{Q}_2^*$ is a determinantal function if 2 is unramified in $E/\mathbf{Q}$. By the calculations of [7], $(\gamma|-) = Ind_{G(E_w/\mathbf{Q}_p)}^{G(E/\mathbf{Q})}(\gamma_2|-)$ and $\gamma_2 = 4(1+\pi)$ where $E_w = \mathbf{Q}_2(\pi)$ so that $\pi = \sqrt{2}$ or $\pi = \sqrt{10}$. Hence

$$(\gamma|\chi) = \begin{cases} 8 & if\ \chi(a_2) = 1, \\ 8\pi & otherwise. \end{cases}$$

and in all cases when 2 ramifies

$$\prod_{\chi=1,\chi_1,\chi_2,\chi_1\chi_2} \frac{f(\chi)^{1/2}}{(\gamma \mid \chi)} = 2^a 5^b$$

for suitable integers $a$ and $b$. Since $2^a 5^b 2^{-v_2(2^a 5^b)} = 5^b \equiv 1$ (modulo 4) the torsion component of $T_S(\mathcal{L},1)$ may also be computed by replacing its Hom-description representative by the simplified function

$$\chi \mapsto (\prod_{p\in S_f} det(1-p^{-1}Frob_w^{-1} \mid \chi^{I(E_w/\mathbf{Q}_p)})^{-1}.$$

Therefore the torsion component of

$$\sum_{p\in S_f,\ G(E_w/\mathbf{Q}_p)=G(E/\mathbf{Q})} \chi_w^\bullet(\mathcal{L},V_w) - T_S(\mathcal{L},1)$$

is equal to

$$\prod_{\substack{p \in S_f, \ G(E_w/\mathbf{Q}_p) \neq G(E/\mathbf{Q}) \\ \chi = 1, \chi_1, \chi_2, \chi_1\chi_2}} det(1 - p^{-1} Frob_w^{-1} \mid \chi^{I(E_w/\mathbf{Q}_p)}) 2^{-v_2(det(1 - p^{-1} Frob_w^{-1} \mid \chi^{I(E_w/\mathbf{Q}_p)}))}$$

modulo 4. However this expression is a square in $(\mathbf{Z}/4)^*$ and is therefore trivial.
□